# Weak stability and generalized weak convolution for random vectors and stochastic processes

**Jolanta K. Misiewicz**[1]

*University of Zielona Góra*

**Abstract:** A random vector $\mathbf{X}$ is weakly stable iff for all $a, b \in \mathbb{R}$ there exists a random variable $\Theta$ such that $a\mathbf{X} + b\mathbf{X}' \stackrel{d}{=} \mathbf{X}\Theta$. This is equivalent (see [11]) with the condition that for all random variables $Q_1, Q_2$ there exists a random variable $\Theta$ such that

$$\mathbf{X}Q_1 + \mathbf{X}'Q_2 \stackrel{d}{=} \mathbf{X}\Theta,$$

where $\mathbf{X}, \mathbf{X}', Q_1, Q_2, \Theta$ are independent. In this paper we define generalized convolution of measures defined by the formula

$$\mathcal{L}(Q_1) \oplus_\mu \mathcal{L}(Q_2) = \mathcal{L}(\Theta),$$

if the equation (∗) holds for $\mathbf{X}, Q_1, Q_2, \Theta$ and $\mu = \mathcal{L}(\Theta)$. We study here basic properties of this convolution, basic properties of $\oplus_\mu$-infinitely divisible distributions, $\oplus_\mu$-stable distributions and give a series of examples.

## 1. Introduction

The investigations of weakly stable random variables started in seventies by the papers of Kucharczak and Urbanik (see [8, 15]). Later there appeared a series of papers on weakly stable distributions written by Urbanik, Kucharczak and Vol'kovich (see e.g. [9, 16, 17, 18]). Recently there appeared a paper written by Misiewicz, Oleszkiewicz and Urbanik (see [11]), where one can find a full characterization of weakly stable distributions with non-trivial discrete part, and a substantial attempt to characterize weakly stable distributions in general case.

In financial mathematics, insurance mathematics and many different areas of science people are trying to predict future behaviour of certain processes by stochastic modelling. Using independent random variables in a variety of constructions turned out to be not sufficient for modelling real events. Multidimensional stable distributions have nice linear properties and enable more complicated structures of dependencies. On the other hand stable distributions are very difficult in calculations because of the complicated density functions and because of the possibility of unbounded jumps stable stochastic processes.

Also the distributions called copulas, extensively investigated recently, are giving the possibility of modelling complicated structure of dependencies. Namely, for every choice of parameters such as covariance matrix or a wines structure of conditional dependency coefficients one can find an arbitrarily nice copula with this parameters.

[1]Department of Mathematics Informatics and Econometry, University of Zielona Góra, ul. Podgórna 50, 65-246 Zielona Góra, Poland, e-mail: j.misiewicz@wmie.uz.zgora.pl

*AMS 2000 subject classifications:* 60A10, 60B05, 60E05, 60E07, 60E10.

*Keywords and phrases:* weakly stable distribution, symmetric stable distribution, $\ell_\alpha$-symmetric distribution, scale mixture.

In this situation weakly stable distributions and processes seem to be good candidates for using in stochastic modelling. They have nice linear properties, i.e. if $(\mathbf{X}_i)$ is a sequence of independent identically distributed random vectors with the weakly stable distribution $\mu$ then every linear combination $\sum a_i \mathbf{X}_i$ has the same distribution as $\mathbf{X}_1 \cdot \Theta$ for some random variable $\Theta$ independent of $\mathbf{X}_1$. This condition holds not only when $(a_i)$ is a sequence of real numbers, but also when $(a_i)$ is a sequence of random variables for $a_i, \mathbf{X}_i$, $i = 1, 2, \dots$, mutually independent. This means that dependence structure of the linear combination $\sum a_i \mathbf{X}_i$ and dependence structure of the random vector $\mathbf{X}_1$ are the same, and the sequence $(a_i)$ is responsible only for the radial behaviour. Moreover weak stability is preserved under taking linear operators $A(\mathbf{X}_1)$, under taking projections or functionals $\langle \xi, \mathbf{X}_1 \rangle$. On the other hand radial properties of a distribution can be arbitrarily defined by choosing a proper random variable $\Theta$ independent of $\mathbf{X}_1$ and considering the distribution of $\Theta \cdot \mathbf{X}_1$. Similar properties of tempered stable distributions (see e.g. [12]) are the reason why they are so important now in statistical physics to model turbulence, or in mathematical finance to model stochastic volatility.

Throughout this paper we denote by $\mathcal{L}(\mathbf{X})$ the distribution of the random vector $\mathbf{X}$. If random vectors $\mathbf{X}$ and $\mathbf{Y}$ have the same distribution we will write $\mathbf{X} \stackrel{d}{=} \mathbf{Y}$. By $\mathcal{P}(\mathbb{E})$ we denote the set of all probability measures on a Banach space (or on a set) $\mathbb{E}$. We will use the simplified notation $\mathcal{P}(\mathbb{R}) = \mathcal{P}$, $\mathcal{P}([0, +\infty)) = \mathcal{P}^+$.

For every $a \in \mathbb{R}$ and every probability measure $\mu$ we define the rescaling operator $T_a : \mathcal{P}(\mathbb{E}) \to \mathcal{P}(\mathbb{E})$ by the formula:

$$T_a \mu(A) = \begin{cases} \mu(A/a) & \text{for } a \neq 0; \\ \delta_0(A) & \text{for } a = 0, \end{cases}$$

for every Borel set $A \in \mathbb{E}$. Equivalently $T_a \mu$ is the distribution of the random vector $a\mathbf{X}$ if $\mu$ is the distribution of the vector $\mathbf{X}$.

The scale mixture $\mu \circ \lambda$ of a measure $\mu \in \mathcal{P}(\mathbb{E})$ with respect to the measure $\lambda \in \mathcal{P}$ is defined by the formula:

$$\mu \circ \lambda(A) \stackrel{def}{=} \int_{\mathbb{R}} T_s \mu(A) \, \lambda(ds).$$

It is easy to see that $\mu \circ \lambda$ is the distribution of random vector $\mathbf{X}\Theta$ if $\mu = \mathcal{L}(\mathbf{X})$, $\lambda = \mathcal{L}(\Theta)$, $\mathbf{X}$ and $\Theta$ are independent. In the language of characteristic functions we obtain

$$\widehat{\mu \circ \lambda}(\mathbf{t}) = \int_{\mathbb{R}} \widehat{\mu}(\mathbf{t}s) \lambda(ds).$$

It is known that for a symmetric random vector $\mathbf{X}$ independent of random variable $\Theta$ we have $\mathbf{X}\Theta \stackrel{d}{=} \mathbf{X}|\Theta|$. From this property we obtain that if $\mu$ is a symmetric probability distribution then

$$\mu \circ \lambda = \mu \circ |\lambda|,$$

where $|\lambda| = \mathcal{L}(|\Theta|)$.

**Definition 1.** A probability measure $\mu \in \mathcal{P}(\mathbb{E})$ is weakly stable (or weakly stable on $[0, \infty)$) if for every choice of $\lambda_1, \lambda_2 \in \mathcal{P}$ ($\lambda_1, \lambda_2 \in \mathcal{P}_+$) there exists $\lambda \in \mathcal{P}$ ($\lambda \in \mathcal{P}_+$) such that

$$(\lambda_1 \circ \mu) * (\lambda_2 \circ \mu) = \lambda \circ \mu.$$

If $\mu$ is not symmetric then the measure $\lambda$ is uniquely determined. This fact was proven in [11] in the case of a weakly stable measure $\mu$, and in [15] in the case of $\mu$ weakly stable on $[0, \infty)$. If the measure $\mu$ is symmetric then only the symmetrization of $\lambda$ is uniquely determined (see [11], Remark 1). In this case we can always replace the measure $\lambda$ by its symmetrization $(\frac{1}{2}\delta_1 + \frac{1}{2}\delta_{-1}) \circ \lambda$. For the convenience in this paper we will assume that for symmetric $\mu$ the measure $\lambda$ is concentrated on $[0, \infty)$ taking if necessary $|\lambda|$ instead of $\lambda$.

The most important Theorem 1 in the paper [11] states that the distribution $\mu$ is a weakly stable if and only if for every $a, b \in \mathbb{R}$ there exists a probability distribution $\lambda \in \mathcal{P}$ such that $T_a\mu * T_b\mu = \mu \circ \lambda$. Moreover we know that (Th. 6 in [11]) if $\mu$ is weakly stable probability measure on a separable Banach space $\mathbb{E}$ then either there exists $a \in \mathbb{E}$ such that $\mu = \delta_a$, or there exists $a \in \mathbb{E} \setminus \{0\}$ such that $\mu = \frac{1}{2}(\delta_a + \delta_{-a})$, or $\mu(\{a\}) = 0$ for every $a \in \mathbb{E}$.

Many interesting classes of weakly stable distributions are already known in the literature. Symmetric stable random vectors are weakly stable, strictly stable vectors are weakly stable on $[0, \infty)$. Uniform distributions on the unit spheres $S^{n-1} \subset \mathbb{R}^n$, their k-dimensional projections and their linear deformations by linear operators are weakly stable. The beautiful class of extreme points in the set of $\ell_1$-symmetric distributions in $\mathbb{R}^n$ given by Cambanis, Keener and Simons (see [4]) is weakly stable.

## 2. Generalized weak convolution

**Definition 2.** Let $\mu \in \mathcal{P}(\mathbb{E})$ be a nontrivial weakly stable measure, and let $\lambda_1, \lambda_2$ be probability measures on $\mathbb{R}$. If

$$(\lambda_1 \circ \mu) * (\lambda_2 \circ \mu) = \lambda \circ \mu,$$

then the generalized convolution of the measures $\lambda_1, \lambda_2$ with respect to the measure $\mu$ (notation $\lambda_1 \oplus_\mu \lambda_2$) is defined as follows

$$\lambda_1 \oplus_\mu \lambda_2 = \begin{cases} \lambda & \text{if } \mu \text{ is not symmetric;} \\ |\lambda| & \text{if } \mu \text{ is symmetric.} \end{cases}$$

If $\Theta_1, \Theta_2$ are random variables with distributions $\lambda_1, \lambda_2$ respectively then the random variable with distribution $\lambda_1 \oplus_\mu \lambda_2$ we will denote by $\Theta_1 \oplus_\mu \Theta_2$. Thus we have

$$\Theta_1 \mathbf{X}' + \Theta_2 \mathbf{X}'' \stackrel{d}{=} (\Theta_1 \oplus_\mu \Theta_2) \mathbf{X},$$

where $\mathbf{X}, \mathbf{X}', \mathbf{X}''$ have distribution $\mu$, $\Theta_1, \Theta_2, \mathbf{X}', \mathbf{X}''$ and $\Theta_1 \oplus_\mu \Theta_2, \mathbf{X}$ are independent. One can always choose such versions of $\Theta_1 \oplus_\mu \Theta_2$ and $\mathbf{X}$ that the above equality holds almost everywhere.

Now it is easy to see that the following lemma holds.

**Lemma 1.** *If the weakly stable measure $\mu \in \mathcal{P}(\mathbb{E})$ is not trivial then*

(1) $\lambda_1 \oplus_\mu \lambda_2$ *is uniquely determined;*
(2) $\lambda_1 \oplus_\mu \lambda_2 = \lambda_2 \oplus_\mu \lambda_1$;
(3) $\lambda \oplus_\mu \delta_0 = \lambda$;
(4) $(\lambda_1 \oplus_\mu \lambda_2) \oplus_\mu \lambda_3 = \lambda_1 \oplus_\mu (\lambda_2 \oplus_\mu \lambda_3)$;
(5) $T_a(\lambda_1 \oplus_\mu \lambda_2) = (T_a\lambda_1) \oplus_\mu (T_a\lambda_2)$.

**Example 1.** It is known that the random vector $\mathbf{U}^n = (U_1, \ldots, U_n)$ with the uniform distribution $\omega_n$ on the unit sphere $S_{n-1} \subset \mathbb{R}^n$ is weakly stable. The easiest way to see this is using the characterizations of a rotationally invariant vectors.

Let us recall that the random vector $\mathbf{X} \in \mathbb{R}^n$ is rotationally invariant (spherically symmetric) if $L(\mathbf{X}) \stackrel{d}{=} \mathbf{X}$ for every unitary linear operator $L : \mathbb{R}^n \to \mathbb{R}^n$. It is known (see [5, 14] for the details) that the following conditions are equivalent

(a) $\mathbf{X} \in \mathbb{R}^n$ is rotationally invariant
(b) $\mathbf{X} \stackrel{d}{=} \Theta \mathbf{U}^n$, where $\Theta = \|\mathbf{X}\|_2$ is independent of $\mathbf{U}^n$,
(c) the characteristic function of $\mathbf{X}$ has the form

$$\mathbf{E} e^{i \langle \xi, \mathbf{X} \rangle} = \varphi_{\mathbf{X}}(\xi) = \varphi(\|\xi\|_2)$$

for some symmetric function $\varphi : \mathbb{R} \to \mathbb{R}$.

Now let $\mathcal{L}(\Theta_1) = \lambda_1$, $\mathcal{L}(\Theta_2) = \lambda_2$ be such that $\Theta_1, \Theta_2, \mathbf{U}^{n1}, \mathbf{U}^{n2}$ are independent, $\mathbf{U}^{n1} \stackrel{d}{=} \mathbf{U}^{n2} \stackrel{d}{=} \mathbf{U}^n$. In order to prove weak stability of $\mathbf{U}^n$ we consider the characteristic function $\psi$ of the vector $\Theta_1 \mathbf{U}^{n1} + \Theta_2 \mathbf{U}^{n2}$

$$\begin{aligned} \psi(\xi) &= \mathbf{E} \exp\{i \langle \xi, \Theta_1 \mathbf{U}^{n1} + \Theta_2 \mathbf{U}^{n2} \rangle\} \\ &= \mathbf{E} \exp\{i \langle \xi, \Theta_1 \mathbf{U}^{n1} \rangle\} \mathbf{E} \exp\{i \langle \xi, \Theta_2 \mathbf{U}^{n2} \rangle\} = \varphi_1(\|\xi\|_2) \varphi_2(\|\xi\|_2). \end{aligned}$$

It follows from the condition (c) that $\Theta_1 \mathbf{U}^{n1} + \Theta_2 \mathbf{U}^{n2}$ is also rotationally invariant. Using condition (b) we obtain that $\Theta_1 \mathbf{U}^{n1} + \Theta_2 \mathbf{U}^{n2} \stackrel{d}{=} \Theta \mathbf{U}^n$ for some random variable $\Theta$, which we denote by $\Theta_1 \oplus_{\omega_n} \Theta_2$. This means that $\mathbf{U}^n$ is weakly stable and

$$\begin{aligned} \Theta_1 \oplus_{\omega_n} \Theta_2 &= \|\Theta_1 \mathbf{U}^{n1} + \Theta_2 \mathbf{U}^{n2}\|_2 \\ &= \left( \sum_{k=1}^n (\Theta_1 U_k^{n1} + \Theta_2 U_k^{n2})^2 \right)^{1/2}, \end{aligned}$$

where $\mathbf{U}^{ni} = (U_1^{ni}, \ldots, U_n^{ni})$, $i = 1, 2$. Since $\mathbf{U}^2 = (\cos \varphi, \sin \varphi)$ for the random variable $\varphi$ with uniform distribution on $[0, 2\pi]$, then in the case $n = 2$ we get

$$\Theta_1 \oplus_{\omega_n} \Theta_2 = \left( \Theta_1^2 + \Theta_2^2 + 2\Theta_1 \Theta_2 \cos(\alpha - \beta) \right)^{1/2},$$

where $\Theta_1, \Theta_2, \alpha, \beta$ are independent, $\alpha$ and $\beta$ have uniform distribution on the interval $[0, 2\pi]$. It is easy to check that $\cos(\alpha - \beta)$ has the same distribution as $\cos(\alpha)$, thus we have

$$\Theta_1 \oplus_{\omega_n} \Theta_2 \stackrel{d}{=} \left( \Theta_1^2 + \Theta_2^2 + 2\Theta_1 \Theta_2 \cos(\alpha) \right)^{1/2}.$$

**Definition 3.** Let $\mathcal{L}(\Theta) = \lambda$, and let $\mu = \mathcal{L}(\mathbf{X})$ be a weakly stable measure on $\mathbb{E}$. We say that the measure $\lambda$ (random variable $\Theta$) is $\mu$-weakly infinitely divisible if for every $n \in \mathbb{N}$ there exists a probability measure $\lambda_n$ such that

$$\lambda = \lambda_n \oplus_\mu \cdots \oplus_\mu \lambda_n, \quad (n\text{-times}),$$

where (for the uniqueness) $\lambda_n \in \mathcal{P}_+$ if $\mu$ is weakly stable on $[0, \infty)$ or if $\mu$ is symmetric, and $\lambda_n \in \mathcal{P}$ if $\mu$ is weakly stable nonsymmetric.

Notice that if $\lambda$ is $\mu$-weakly infinitely divisible then the measure $\lambda \circ \mu$ is infinitely divisible in the usual sense. However if $\lambda \circ \mu$ is infinitely divisible then it does not have to imply $\mu$-infinite divisibility of $\lambda$.

**Example 2.** If $\gamma_\alpha$ is a strictly $\alpha$-stable (symmetric $\alpha$-stable) distribution on a separable Banach space $\mathbb{E}$ then it is weakly stable on $[0,\infty)$ (weakly stable). Simple application of the definition of stable distribution shows that

$$\Theta_1 \oplus_{\gamma_\alpha} \Theta_2 \stackrel{d}{=} (\Theta_1^\alpha + \Theta_2^\alpha)^{1/\alpha}, \quad \left(\Theta_1 \oplus_{\gamma_\alpha} \Theta_2 \stackrel{d}{=} (|\Theta_1|^\alpha + |\Theta_2|^\alpha)^{1/\alpha}\right).$$

Now we see that $\Theta$ is $\gamma_\alpha$- weakly infinitely divisible if and only if $\Theta^\alpha$ (respectively $|\Theta|^\alpha$) is infinitely divisible in the usual sense.

**Lemma 2.** *Let $\mu$ be a weakly stable distribution, $\mu \neq \delta_0$. If $\lambda$ is $\mu$-weakly infinitely divisible then there exists a family $\{\lambda^r : r \geq 0\}$ such that*
(1) $\lambda^0 = \delta_0$, $\lambda^1 = \lambda$;
(2) $\lambda^r \oplus_\mu \lambda^s = \lambda^{r+s}$, $r, s \geq 0$;
(3) $\lambda^r \Rightarrow \delta_0$ *if* $r \to 0$.

*Proof.* If $\lambda$ is $\mu$-weakly infinitely divisible then for every $n \in \mathbb{N}$ there exists a measure $\lambda_n$ such that

$$(\lambda_n \circ \mu)^{*n} = \lambda \circ \mu,$$

where $\nu^{*n}$ denotes the $n$'th convolution power of the measure $\nu$. We define $\lambda^{1/n} := \lambda_n$. Weak stability of the measure $\mu$ implies that for every $k, n \in \mathbb{N}$ there exists a probability measure which we denote by $\lambda^{k/n}$ such that

$$\lambda^{k/n} \circ \mu = \left(\lambda^{1/n} \circ \mu\right)^{*k} = (\lambda \circ \mu)^{*k/n}.$$

The last expression follows from the infinite divisibility of the measure $\lambda \circ \mu$. We see here that for every $n, k, m \in \mathbb{N}$ we have

$$\lambda^{km/nm} = \lambda^{k/n},$$

since

$$(\lambda \circ \mu)^{*km/nm} = (\lambda \circ \mu)^{*k/n}.$$

Now let $x > 0$ and let $(r_n)_n$ be a sequence of rational numbers such that $r_n \to x$ when $n \to \infty$. Since $(\lambda \circ \mu)^{*r_n} \to (\lambda \circ \mu)^{*x}$ and

$$\{\lambda^{r_n} \circ \mu : n \in \mathbb{N}\} = \{(\lambda \circ \mu)^{*r_n} : n \in \mathbb{N}\}$$

then this family of measures is tight. Lemma 2 in [11] implies that also the family $\{\lambda^{r_n} : n \in \mathbb{N}\}$ is tight, so there exists a subsequence $\lambda^{r_{n_k}}$ weakly convergent to a probability measure which we call $\lambda^x$. Since $\lambda^x \circ \mu = (\lambda \circ \mu)^{*x}$ then uniqueness of the measure $\lambda^x$ follows from the uniqueness of $(\lambda \circ \mu)^{*x}$, Remark 1 in [11] and our assumptions.

To see (3) let $r_n \to 0$, $r_n > 0$. Since

$$\lambda^{r_n} \circ \mu = (\lambda \circ \mu)^{r_n} \Rightarrow \delta_0 = \delta_0 \circ \mu,$$

then $\{(\lambda \circ \mu)^{r_n} : n \in \mathbb{N}\}$ is tight, and by Lemma 2 in [11] the set $\{\lambda^{r_n} : n \in \mathbb{N}\}$ is also tight. Let $\{r'_n\}$ be the subsequence of $\{r_n\}$ such that $\lambda^{r'_n}$ converges weakly to some probability measure $\lambda^0$. Then we have

$$\lambda^{r'_n} \circ \mu \Rightarrow \lambda^0 \circ \mu,$$

and therefore $\lambda^0 \circ \mu = \delta_0 \circ \mu$. If $\mu$ is not symmetric then Remark 1 in [11] implies that $\lambda^0 = \delta_0$. If $\mu$ is symmetric then by our assumptions $\lambda$ and $\lambda^{r_n}$ are concentrated on $[0, \infty)$, thus also $\lambda^0$ is concentrated on $[0, \infty)$. Since by Remark 1 in [11] the symmetrization of the mixing measure is uniquely determined in this case we also conclude that $\lambda^0 = \delta_0$. □

## 3. $\mu$-weakly stable random variables and vectors

**Definition 4.** Let $\mu$ be a weakly stable distribution, $\mu \neq \delta_0$. We say that the probability measure $\lambda$ is $\mu$-weakly stable if

$$\forall\, a, b > 0 \,\exists\, c > 0 \ \text{ such that } \ T_a \lambda \oplus_\mu T_b \lambda = T_c \lambda.$$

We say that the random variable $\Theta$ is $\mu$-weakly stable if

$$\forall\, a, b > 0 \,\exists\, c > 0 \ \text{ such that } \ (a\Theta)\mathbf{X} + (b\Theta')\mathbf{X}' \stackrel{d}{=} (c\Theta)\mathbf{X},$$

where the random variable $\Theta'$ is an independent copy of $\Theta$, the vectors $\mathbf{X}$ and $\mathbf{X}'$ have distribution $\mu$ and $\Theta, \Theta', \mathbf{X}, \mathbf{X}'$ are independent.

Directly from the definition we see the following

**Lemma 3.** *Let $\mu$ be a weakly stable distribution, $\mu \neq \delta_0$. A probability measure $\lambda$ is $\mu$-weakly stable iff the measure $\lambda \circ \mu$ is strictly stable in the usual sense, thus there exists $\alpha \in (0, 2]$ such that $\lambda \circ \mu$ is strictly $\alpha$-stable. In such a case we will say that $\lambda$ is $\mu$-weakly $\alpha$-stable.*

*Proof.* Let us define

$$\psi(t) = \mathbf{E} e^{it\mathbf{X}\Theta} = \int \widehat{\mu}(ts) \lambda(ds),$$

where $\mathbf{X}$ has distribution $\mu$, $\Theta$ has distribution $\lambda$, $\mathbf{X}$ and $\Theta$ are independent. The condition of $\mu$-weak stability of $\lambda$ can be written in the following way

$$\forall\, a, b > 0 \,\exists\, c > 0 \ \text{ such that } \ \psi(at)\psi(bt) = \psi(ct),$$

which is the functional equation defining strictly stable characteristic functions. $\square$

**Example 3.** Let $\mathbf{U}^n$ be a random vector with the uniform distribution $\mu = \omega_n$ on the unit sphere $S_{n-1} \subset \mathbb{R}^n$ (or $\mathbf{U}^{n,k}$-any its projection into $\mathbb{R}^k$, $k < n$). The $\omega_n$-weakly Gaussian random variable $\Gamma_n$ is defined by the following equation:

$$\mathbf{U}^n \cdot \Gamma_n \stackrel{d}{=} (X_1, \ldots, X_n) = \mathbf{X},$$

where $\mathbf{U}^n$ and $\Gamma_n$ are independent, $\mathbf{X}$ is an $n$-dimensional Gaussian random vector with independent identically distributed coordinates. For convenience we can assume that each $X_i$ has distribution $N(0, 1)$. It follows from the condition (b) in the characterization of rotationally invariant random vectors given in Example 1 that $\Gamma_n$ has the same distribution as $\|\mathbf{X}\|_2$. Simple calculations show that $\Gamma_n$ has density

$$f_{2,n}(r) = \frac{2}{2^{n/2} \Gamma(\frac{n}{2})} r^{n-1} e^{-r^2/2}.$$

For $n = 2$ this is a Rayleigh distribution with parameter $\lambda = 2$, thus the Rayleigh distribution is $\omega_2$-weakly Gaussian. For $n = 3$ this is the Maxwell distribution with parameter $\lambda = 2$, thus Maxwell distribution is $\omega_3$-weakly Gaussian. Let us remind that the generalized Gamma distribution with parameters $\lambda, p, a > 0$ (notation $\Gamma(\lambda, p, a)$) is defined by its density function

$$f(x) = \frac{a}{\Gamma(p/a) \lambda^{p/a}} x^{p-1} \exp\left\{-\frac{x^a}{\lambda}\right\}, \quad \text{for } x > 0.$$

Thus we have that the generalized Gamma distribution $\Gamma(\lambda, n, 2)$ is $\omega_n$-weakly Gaussian.

Now let $\Theta_\alpha^n$ be an $\omega_n$-weakly $\alpha$-stable random variable. Then $\mathbf{U}^n \Theta_\alpha^n$ is a rotationally invariant $\alpha$-stable random vector for $\mathbf{U}^n$ independent of $\Theta_\alpha^n$. On the other hand, every rotationally invariant $\alpha$-stable random vector has the same distribution as $\mathbf{Y}\sqrt{\Theta_{\alpha/2}}$, where $\mathbf{Y}$ is rotationally invariant Gaussian random vector independent of the nonnegative variable $\Theta_{\alpha/2}$ with the Laplace transform $e^{-t^{\alpha/2}}$. Finally we have

$$\mathbf{U}^n \cdot \Theta_\alpha^n \stackrel{d}{=} \mathbf{U}^n \Gamma_n \sqrt{\Theta_{\alpha/2}},$$

for $\mathbf{U}^n$, $\Gamma_n$ and $\Theta_{\alpha/2}$ independent. This implies that the density of a $\omega_n$-weakly $\alpha$-stable random variable $\Theta_\alpha^n$ is given by

$$f_{\alpha,n}(r) = \int_0^\infty f_{2,n}\left(\frac{r}{\sqrt{s}}\right) \frac{1}{\sqrt{s}} f_{\alpha/2}(s) ds.$$

In particular if we take $\alpha = 1$ then

$$f_{1/2}(s) = \frac{1}{\sqrt{2\pi}} x^{-3/2} e^{-1/(2x)}, \quad x > 0.$$

Simple calculations and the equality $\sqrt{\pi}\Gamma(2s) = 2^{2s-1}\Gamma(s)\Gamma(s+\frac{1}{2})$, $s > 0$, show that

$$f_{1,n}(r) = \frac{2^{2-n}\Gamma(n)}{\Gamma(n/2)\Gamma(n/2)} \frac{r^{n-1}}{(r^2+1)^{(n+1)/2}}, \quad r > 0,$$

is the density function of a $\omega_n$-weakly Cauchy distribution.

**Example 4.** We know that for every symmetric $\alpha$-stable random vector $\mathbf{X}$ with distribution $\gamma_\alpha$ on any separable Banach space $\mathbb{E}$ and every $p \in (0, 1)$ the random vector $\mathbf{X}\Theta_p^{1/\alpha}$ is symmetric $\alpha p$-stable for $\mathbf{X}$ independent of $\Theta_p$ with the distribution $\lambda_p$ and the Laplace transform $e^{-t^p}$. Since symmetric stable vectors are weakly stable we obtain that

$$\forall \gamma_\alpha \ \forall \ p \in (0,1) \quad \lambda_p \text{ is } \gamma_\alpha - \text{weakly } \alpha p - \text{stable}.$$

## 4. $(\lambda, \mu)$-weakly stable Lévy processes

In this section we construct a Lévy process based on a nontrivial weakly stable probability measure $\mu$. The measure $\mu$ can be defined on the real line, on $\mathbb{R}^n$ or on a Banach space $\mathbb{E}$. By $\lambda$ we will denote in this section a $\mu$-weakly infinitely divisible distribution on $\mathbb{R}$.

Let $T = [0, \infty)$ and let $m$ be a Borel measure on $T$. We say that $\{\mathbf{X}_t : t \in T\}$ is a $(\lambda, \mu)$-weakly stable Lévy process if the following conditions hold:

(a) $\mathbf{X}_0 \equiv 0$;
(b) $\mathbf{X}_t$ has independent increments;
(c) $\mathbf{X}_t$ has distribution $\lambda^{m[0,t)} \circ \mu$.

If $m$ is equal to the Lebesgue measure on $T$ then this process has stationary increments.

Notice that for $\lambda = \delta_1$ and $\mu = N(0, 1)$ we obtain with this construction the Brownian motion.

**Example 5.** Let $\gamma_\alpha$ be a strictly stable distribution on a separable Banach space $\mathbb{E}$ with the characteristic function $\exp\{-R(\xi)\}$ and let $\Theta_r$, $r > 0$, denote the random variable with distribution

$$\mathbf{P}\{\Theta_r = k\} = \binom{r+k-1}{k}(1-p)^k p^r, \quad k = 0, 1, 2 \ldots,$$

for some $p \in (0,1)$. Since $R(t\xi) = |t|^\alpha R(\xi)$ it is easy to see that the measure $\lambda = \mathcal{L}(\Theta_1^{1/\alpha})$ is $\gamma_\alpha$-weakly infinitely divisible, $\lambda^r = \mathcal{L}(\Theta_r^{1/\alpha})$, and

$$\lambda^r \oplus_{\gamma_\alpha} \lambda^s = \lambda^{r+s}.$$

Let $\{\mathbf{X}_t : t \in T\}$ be the $(\lambda, \gamma_\alpha)$-weakly stable Lévy process. Then $\mathbf{X}_t$ has the following characteristic function:

$$\begin{aligned}
\mathbf{E}\exp\{i\langle \xi, \mathbf{X}_t\rangle\} &= \mathbf{E}\exp\{-R(\xi)\Theta_{m[0,t)}\} \\
&= \sum_{k=0}^\infty \exp\{-R(\xi)k\}\binom{m[0,t)+k-1}{k}(1-p)^k p^{m[0,t)} \\
&= \left(\frac{p}{1-(1-p)\exp\{-R(\xi)\}}\right)^{m[0,t)}.
\end{aligned}$$

**Example 6.** For $\gamma_\alpha$ being a strictly stable distribution on a separable Banach space $\mathbb{E}$ with the characteristic function $\exp\{-R(\xi)\}$ and $\lambda = \mathcal{L}(Q^{1/\alpha})$, where $Q$ has Gamma distribution with parameters $p = 1$ and $a > 0$ we obtain $(\lambda, \gamma_\alpha)$-weakly stable Levy process with the distribution defined by the following characteristic function

$$\begin{aligned}
\mathbf{E}\exp\{i\langle \xi, \mathbf{X}_t\rangle\} &= \mathbf{E}\exp\{-R(\xi)Q_{m[0,t)}\} \\
&= \int_0^\infty \exp\{-R(\xi)s\}\frac{a^{m[0,t)}}{\Gamma(m[0,t))}s^{m[0,t)-1}e^{-as}ds \\
&= \left(\frac{a}{a+R(\xi)}\right)^{m[0,t)}.
\end{aligned}$$

To see this it is enough to notice that $\lambda^r = \mathcal{L}(Q_r^{1/\alpha})$, where $Q_r$ has gamma distribution $\Gamma(r, a)$.

## 5. $\mu$-weakly one-dependent processes

Let us recall that the stochastic process $\{\mathbf{Y}_n : n \in \mathbb{N}\}$ taking values in a separable Banach space $\mathbb{E}$ is one-dependent if for each $n \in \mathbb{N}$ the sequences $\{\mathbf{Y}_1, \ldots, \mathbf{Y}_{n-1}\}$ and $\{\mathbf{Y}_{n+1}, \mathbf{Y}_{n+2}, \ldots\}$ are independent. It is evident that if $f \colon \mathbb{R} \mapsto \mathbb{R}$ is a measurable function and $\{\mathbf{Y}_1, \mathbf{Y}_2, \ldots\}$ is a one-dependent process then also the process $\{f(\mathbf{Y}_1), f(\mathbf{Y}_2), \ldots\}$ is one-dependent. A simples possible one-dependent process can be obtained as $\{f(X_i, X_{i+1}) \colon i = 1, 2, \ldots\}$, where $\{X_i\}$ is a sequence of independent (often identically distributed) random variables. A nice counterexample that not all one-dependent processes have this construction is given in [1, 2].

There are several possibilities for constructing one-dependent processes with distributions which are mixtures of a fixed weakly stable measure $\mu$ on a separable Banach space $\mathbb{E}$. In the first of the following examples we give this construction assuming that the mixing measure is $\mu$-weakly infinitely divisible. In the second example this assumption is omitted, but the weakly stable measure $\mu$ must be stable.

**Example 7.** Let $\mu$ be a weakly stable distribution on a separable Banach space $\mathbb{E}$ and let $\lambda$ be a $\mu$-weakly infinitely divisible measure on $\mathbb{R}$. Assume that $m$ is a $\sigma$-finite measure on a rich enough measure space $(S, \mathcal{B})$ such that $S = \bigcup_{n=1}^{\infty} A_n$, $m(A_n) < \infty$, $n = 1, 2, \ldots$ and $A_i \cap A_j = \emptyset$ for $i \neq j$. With each set $A_n$ we connect the random variable $Z_n$ with distribution $\lambda^{m(A_n)}$ such that $Z_1, Z_2, \ldots$ are independent. Let also $\mathbf{X}_1, \mathbf{X}_2, \ldots$ be the sequence of independent identically distributed random vectors with distribution $\mu$. Now we define

$$\mathbf{Y}_n = Z_n \mathbf{X}_n + Z_{n+1} \mathbf{X}_{n+1}, \quad n = 1, 2, \ldots.$$

It is easy to see that $\{\mathbf{Y}_n \colon n \in \mathbb{N}\}$ is a one-dependent process with the distribution

$$\mathcal{L}(\mathbf{Y}_n) = \lambda^{m(A_n \cup A_{n+1})} \circ \mu.$$

This process is stationary if $m(A_i) = m(A_j)$ for all $i, j \in \mathbb{N}$. If $\mu = \omega_k$ then $\{\mathbf{Y}_n \colon n \in \mathbb{N}\}$ is elliptically contoured. If $\mu = \gamma_\alpha$ then $\{\mathbf{Y}_n \colon n \in \mathbb{N}\}$ is $\alpha$-substable. If $\mu = \gamma_\alpha$ and for some $\beta \in (0, 1)$ $\lambda = \mathcal{L}(\Theta_\beta)$, where $\Theta_\beta$ is nonnegative $\beta$-stable random variable, then $\{\mathbf{Y}_n \colon n \in \mathbb{N}\}$ is $\alpha\beta$-stable $\alpha$-substable.

**Example 8.** Assume that $\gamma_p$ is a symmetric $p$-stable distribution on a separable Banach space $\mathbb{E}$, and let $\{Z_n \colon n \in \mathbb{N}\}$ be any one-dependent stochastic process taking values in $[0, \infty)$. For the sequence $\mathbf{X}_1, \mathbf{X}_2, \ldots$ of i.i.d. random vectors with distribution $\mu$ we define
$$\mathbf{Y}_n = \mathbf{X}_n Z_n^{1/p}, \quad n \in \mathbb{N}.$$

Directly from the construction, it follows that the process $\{\mathbf{Y}_n \colon n \in \mathbb{N}\}$ is one-dependent. This process is also $p$-substable and the characteristic function of the linear combination of its components $\sum a_n \mathbf{Y}_n$ can be easily calculated using the Laplace transform for $\{Z_n \colon n \in \mathbb{N}\}$. Namely for every $\xi \in \mathbb{E}^*$ we have

$$\mathbf{E} \exp \left\{ i \langle \xi, \sum a_n \mathbf{Y}_n \rangle \right\} = \mathbf{E} \exp \left\{ i \sum a_n \langle \xi, \mathbf{X}_n \rangle Z_n^{1/p} \right\}$$
$$= \mathbf{E} \exp \left\{ -\sum |a_n|^p \|\Re(\xi)\|_p^p Z_n > \right\},$$

where $\Re$ is the linear operator fom $\mathbb{E}^{ast}$ into some $L_p$-space such that

$$\mathbf{E} \exp \left\{ i \langle \xi, \mathbf{X}_1 \rangle \right\} = \mathbf{E} \exp \left\{ -\|\Re(\xi)\|_p^p \right\}.$$